\RequirePackage[l2tabu, orthodox]{nag}
\documentclass[a4paper]{amsart}
\usepackage{mathrsfs}
\usepackage{amsfonts}
\usepackage{amsmath}
\usepackage{amssymb}
\usepackage{amsthm}
\usepackage{enumerate}
\usepackage{color}
\usepackage{geometry}
\usepackage[pdftex,linkcolor=blue,citecolor=blue,backref=page]{hyperref}
\allowdisplaybreaks
\newtheorem{theorem}{Theorem}[section]
\newtheorem{definition}[theorem]{Definition}

\newtheorem{lemma}[theorem]{Lemma}
\newtheorem{corollary}[theorem]{Corollary}
\newtheorem{proposition}[theorem]{Proposition}
\newtheorem{example}[theorem]{Example}
\newtheorem{remark}[theorem]{Remark}
\newtheorem{question}[theorem]{Question}
\newtheorem{conjecture}[theorem]{Conjecture}

\newcommand{\GK}{\operatorname{GKdim}}

\newcommand{\id}{\operatorname{id}}
\newcommand{\Id}{\operatorname{Id}}

\newcommand{\Irr}{\operatorname{Irr}}

\newcommand{\NN}{\mathbb{N}}

\newcommand{\ZZ}{\mathbb{Z}}

\newcommand{\FF}{\mathbb{F}}
\newcommand{\G}{\mathcal{G}}

\newcommand{\gsb}{Gr\"obner-Shirshov basis}
\newcommand{\gsbs}{Gr\"obner-Shirshov bases}

\newcommand{\fd}{finite dimensional}

\newcommand{\gp}[2]{\operatorname{gp}\langle #1\mid #2 \rangle}
\newcommand{\sgp}[2]{\operatorname{sgp}\langle #1\mid #2\rangle}

\newcommand{\wt}{\operatorname{wt} }

\title[Growth of monomial algebras and Manturov groups]{Growth of associated monomial algebras with application to Manturov groups}

\author[Zhao]{Xiangui Zhao}
\address{School of Mathematics and Statistics,
Huizhou University, Huizhou, Guangdong 516007, China}
\email{zhaoxg@hzu.edu.cn}
\date{\today} 


\begin{document}
\maketitle
\begin{abstract}
  It is well-known that an associative algebra shares the same growth and Gelfand-Kirillov dimension (GK-dimension) as its associated monomial algebra with respect to a degree-lexicographic order.
This article mainly investigates the relationship between the  GK-dimension of an algebra and that of its associated monomial algebra with respect to a monomial order.
We obtain sufficient conditions on a monomial order such that these two algebras have the same GK-dimension.
Our result generalizes the well-known result and has several applications.
In particular, as an application,
we study the growth of Manturov $(k,n)$-groups for positive integers $n>k$.
It is shown that
the  Manturov $(1,n)$-group has growth equal to $0$ for all $n>1$;
the  Manturov $(2,3)$-group has growth equal to $2$;
and, for all $n>k\geq3$, the  Manturov $(k,n)$-group contains a free subgroup of rank $2$
and thus has exponential growth.
\end{abstract}

\section{Introduction and main results}
Let $\FF$ be a field.
Throughout this article, all algebras are assumed to be associative and unital over $\FF$, unless explicitly stated otherwise.

The notion of growth,
introduced by Gelfand and Kirillov \cite{gelfand1966corps}
for algebras and by Milnor \cite{milner1968curvature} for groups,
is a very useful and powerful tool for investigating algebras, modules, groups, operads, and so on.
The Gelfand-Kirillov dimension (GK-dimension, for short) is a numerical characteristic of growth.
The GK-dimension of an algebra $A$ is defined as
\[
\GK(A):=\sup_V{\limsup_{n\to\infty}}\log_n \dim\sum_{i=0}^nV^i
\]
where the supremum is taken over all finite dimensional subspaces $V$ of $A$.
Note that $\GK(A)$ is an invariant of the algebra $A$
but not an invariant of the ring $A$ (see \cite{Wu1991} for details).
We refer the interested readers to \cite{Krause-Lenagan_2000} for basic properties and applications of GK-dimension of algebras, modules and groups.
For recent studies on GK-dimension (and its applications) of other algebraic structures, see, e.g.,
\cite{Bao2020,Bai2021,TZZ21,BXX22,QXZZ23,BCLS24,LMZZ25}.

Suppose $A=\FF\langle X\rangle/I$,
where $\FF\langle X\rangle$ is the free associative algebra generated by a finite set $X$ and $I$ is an ideal of $\FF\langle X\rangle$.
Assume $<$ is a \emph{monomial order} (Definition \ref{def-MonomialOder}) on $X^*$, the free monoid on $X$.
For $f\in \FF\langle X\rangle$, we define $\overline{f}$,
the \emph{leading monomial} of $f$, to be the maximum monomial (with nonzero coefficient) in the expression for $f$ with respect to $<$.
We associate a monomial ideal $J$ to $I$
by letting $J$ be the ideal of $\FF\langle X\rangle$ generated by $\overline{I}:=\{\overline{f}\mid f\in I\}$.
Then we call $\widetilde{A}=\widetilde{A_<}=\FF\langle X\rangle/J$
the \emph{associated monomial algebra} of $A$ with respect to $<$.
It is a standard result that
$\widetilde{A}$ and $A$ have identical growth function (and thus identical GK-dimension)
with respect to a \emph{deg-lex order} on $X^*$ (defined in Example \ref{exam-deg-lex}),
see, e.g., \cite[Remark 4.1]{bell2015growth}, 
or \cite[page 157]{Krause-Lenagan_2000} for the finitely presented case.

The algebra  $\widetilde{A}$ has some advantages over $A$,
for instance, $\widetilde{A}$ provides a more concrete framework for analyzing its growth.
Note that a deg-lex order is a special example of {monomial orders}.
Then the following question is natural.

\begin{question}
For an arbitrary monomial order $<$ on $X^*$,
does the associated monomial algebra $\widetilde{A_<}$ 
have the same growth and GK-dimension as $A$?
\end{question}

Example \ref{exam-GKneq} shows that the answer to the above question is ``no'' in general.
We investigate sufficient conditions on the monomial order such that $\widetilde{A}$ and $A$ have the same growth.
One of our main results can be stated as follows.

\begin{theorem}\label{thm-main-Intro}
  Let $A$ be an algebra generated by a finite set $X$,
  and let $<$ be a monomial order on $X^*$.
  Suppose that there exists a real polynomial $f(x)$ of degree $d=\deg(f)\geq 1$ such that
  for all $ u,v\in X^*$ with $|u|\gg 0$,
  \begin{align*}
   u<v\implies |u|\leq f(|v|),
  \end{align*}
  where $|u|$ denotes the length of $u$.
   Then $$
  \GK(\widetilde{A_<})\leq \GK(A)\leq d\GK(\widetilde{A_<}).
  $$
\end{theorem}

Theorem \ref{thm-main-Intro} has two immediate applications.
First,
if $<$ is a \emph{weighted deg-lex order} (defined in Section \ref{sec-GK-MonomialAlg}),
then
it follows from Theorem \ref{thm-main-Intro} that,
$\GK(A)=\GK(\widetilde{A})$ (Corollary \ref{cor-deg-lex}),
which generalizes the standard result mentioned above.
Second,
Theorem \ref{thm-main-Intro} also implies a dichotomy of the GK-dimension of Ore extensions over the polynomial algebra in one variable
(Corollary \ref{cor-OreExt}).
This dichotomy is an analogue of the dichotomies given in \cite[Theorem 1.2]{zhao2022growth}
for the GK-dimension of generalized Weyl algebras over the polynomial algebra in two variables,
and in \cite{Gupta2017} for the GK-dimension of simple modules over simple differential rings.
However, the dichotomy given by Corollary \ref{cor-OreExt}
cannot be extended to Ore
extensions over polynomial algebra in more than one variable (see Example \ref{exam-OreExt2v}). 

We also apply Theorem \ref{thm-main-Intro} to the study of the growth of 
Manturov $(k,n)$-group $G_n^k$ (defined in Section \ref{sec-GK-MonomialAlg}, also known as the $n$-strand $k$-free braid group) 
for positive integers $n>k$.
Manturov $(k,n)$-groups
were introduced by Manturov \cite{manturov2015non,manturov2015on} in 2015,
which appear naturally as groups describing dynamical systems of $n$ particles.
It turns out that Manturov groups are closely related to braid groups, geometry, topology, and dynamical systems \cite{manturov2015on,kim2016on,manturov2017the}.
We refer interested readers to the survey \cite{MFKN21} for more background and applications of $G^k_n$,
and to \cite{AK20,AlHussein2021,KM25} for some recent results on $G^k_n$.

Recall from \cite[page 140]{Krause-Lenagan_2000} that a group $G$ has \emph{growth equal to $d$} for a real number $d$ if 
the group algebra $\FF G$ has GK-dimension $d$.
Our main result on the growth of Manturov groups is as follows.

\begin{theorem}
  \label{thm-main2}
  Let $G^k_n$ be the Manturov $(k,n)$-group for integers $n>k\geq1$.
  Then the following statements hold.
  \begin{enumerate}[(i)]
    \item For all $n>1$, the group $G^1_n$ has growth equal to $0$.
    \item The group $G^2_3$ has growth equal to $2$,
    and is a finite extension of $\ZZ^2$.
    \item For all $n>k\geq 3$, the group $G^k_n$ contains a free subgroup of rank $2$,
    and consequently has exponential growth.
  \end{enumerate}  
\end{theorem}

Our analysis of the growth of $G^2_3$ employs a specially constructed \emph{infinite} \gsb\ 
of $\FF G^2_3$ with respect to a deg-lex order.
Note that Bai and Chen \cite[Theorem 3]{BC17} (see also \cite[Theorem 1]{AlHussein2021}) obtained a \emph{finite} \gsb\ for $\FF G^2_3$ with respect to an order that does not satisfy the assumption in 
Theorem \ref{thm-main-Intro}.
As a result, the \gsb\ given by Bai and Chen cannot be used directly 
to calculate the GK-dimension of $\FF G^2_3$ (see Remark \ref{remark-tower}).
As pointed out in \cite{manturov2017the},
the word problem of $G_n^k$ is ``extremely important''.
The word problems are still open for all $G_n^k$ except $G_n^2$, $G_4^3$ and $G_5^4$ (see \cite{FKM20}).
In the case of $G_n^2$, $G_4^3$ and $G_5^4$,
the word problems were solved either geometrically or via connections with other groups (see \cite{manturov2017the,FKM20}).
In this paper, we actually give a pure algebraic way to solve the word problems for $G_3^2$ by using the specially constructed infinite \gsb\ of $\FF G^2_3$
(Corollary \ref{cor-WordProb23}).
It is worth mentioning that
\gsb\ theory is a powerful tool for solving word problems of groups, semigroups, associative algebras, Lie algebras, etc.,
and there have been lots of papers in this direction, see e.g. \cite{shirshov1962Lie,bokut2007markov,chen2009a,BOKUT_2009,BCZ24}.

At present we do not know the growth of $G^2_n$ for $n\geq4$ (see Conjecture \ref{conj}).

The rest of the paper is organized as follows.
In Section \ref{sec-pre},
we introduce basic notions and facts on \gsbs.
Then we study the growth of associated monomial algebras and prove Theorem \ref{thm-main-Intro} in Section \ref{sec-GK-MonomialAlg}.
Finally we apply these algebraic techniques to 
investigate the growth of Manturov groups and prove Theorem \ref{thm-main2} in Section \ref{sec-Manturov}.
\section{Preliminaries}
\label{sec-pre}
In this section, we briefly introduce Gr\"{o}bner-Shirshov bases for associative algebras and semigroups.
For more details we refer to \cite{shirshov1962Lie,bokut1976embeddings},
or to the survey \cite{bokut2014grobner}.

Let $\FF$ be a field,
$\FF\langle X\rangle$ be the free associative algebra
over $\FF$ generated by $X$,
and $ X^{*}$ be the free monoid generated by $X$,
where the empty word is the identity $1$.
An element in $X^*$ is called a \emph{word} or a \emph{monomial},
while an element in $\FF\langle X\rangle$ is called a (noncommutative) \emph{polynomial}.
We denote the \emph{length} (or \emph{degree}) of a word $w$ by $|w|$.
For a subset $S\subseteq X^*$,
denote $S^{\leq n}:=\{u\in S\mid |u|\leq n\}$.

We always suppose that the set $X$ is well-ordered.

Let $X^*$ be a well-ordered set and $f\in \FF\langle X\rangle$. 
Denote by $\overline{f}$ the {\it leading monomial} of $f$,
that is, the maximum monomial (with nonzero coefficient) in $f$. 
For a subset $R\subseteq \FF\langle X\rangle$, 
denote $\overline{R}:=\{\overline{r}\mid r\in R\}$.
We call $f$ {\it monic} if
$\overline{f}$ has coefficient 1,
and call $R\subseteq \FF\langle X\rangle$ \emph{monic} if all elements in $R$ are monic.
Denote by $\FF R$ the subspace of $\FF\langle X\rangle$ spanned by $R$,
and by $\Id(R)$ the ideal of $\FF\langle X\rangle$ generated by $R$.

\begin{definition}\label{def-MonomialOder}
  A \emph{monomial order} $<$ on $X^*$ is a well order on $X^*$ such that, 
for $u, v\in X^*$, 
$$
u < v \implies w_{1}uw_{2} < w_{1}vw_{2},  \ \mbox{for all }
 w_{1}, w_{2}\in  X^*.
$$
\end{definition}

\begin{example}\label{exam-deg-lex}
A standard example of monomial orders on $X^*$ is the  degree-lexicographic (deg-lex, for short) order, which compares two words first by degree (length) and then lexicographically.
\end{example}

Another example of monomial orders is a weighted deg-lex order (defined in Section \ref{sec-GK-MonomialAlg}),
which is a generalization of a deg-lex order.

\begin{example}\label{exam-Tower}
  For $1\neq u\in X^*$, let $u':=\max u$ be maximum letter appearing in $u$.
  Then $u$ has a unique decomposition 
  $$
  u=u_0u'u_1u'u_2\cdots u_{n_u-1}u'u_{n_u},
  $$
  where $u_i\in (X\setminus \{u'\})^*$ for $0\leq i\leq n_u$ and ${n_u}\in \NN^*$.
  Let $\wt(u)=(u',n_u,u_0,u_1,\dots,u_n)$.
  For $ 1\neq v\in X^*$,
  set $1<v$, and set $u<v$ if $\wt(u)<\wt(v)$ lexicographically,
  where $u_i<v_i$ for $0\leq i\leq u_n$ if and only if 
  $$
  u_i=1<v_i, \text{ or }
  u_i,v_i\neq 1 \text{ and }\wt(u_i)<\wt(v_i) \text{ recursively}.
  $$
  This order is called a \emph{tower order} on $X^*$.
  It is routine to check that a tower order is a monomial order.
\end{example}

Let $<$ be a monomial order on $X^*$. 
Suppose  $f$ and $g$ are two monic polynomials in $\FF\langle X\rangle$.
We define two kinds of {\it compositions} of $f$ and $g$ as follows.

\begin{enumerate}[(i)]
  \item
  If $w$ is a word such that $w=\overline{f}b=a\overline{g}$ for some
$a,b\in X^*$ with $|\overline{f}|+|\overline{g}|>|w|$, then the polynomial
 $ (f,g)_w:=fb-ag$ is called the {\it intersection composition} of $f$ and
$g$ with respect to $w$.
\item
If  $w=\overline{f}=a\overline{g}b$ for some $a,b\in X^*$, then the
polynomial $ (f,g)_w:=f - agb$ is called the {\it inclusion
composition} of $f$ and $g$ with respect to $w$.
\end{enumerate}

The word $w$ is called the \emph{ambiguity} of the composition
$(f,g)_w$ in each case.

Let $S\subseteq \FF\langle X\rangle$ be monic.
A polynomial $p\in \FF\langle X\rangle$ is called {\it trivial} modulo $(S,w)$ and write
\(
p\equiv0  \mod (S,w)
\)
if $p$ can be written as a finite sum
\[
p=\sum\alpha_i a_i s_i b_i,\ \text{each } \alpha_i\in \FF,
 a_i,b_i\in X^{*}, \ s_i\in S, \ a_i \overline{s_i} b_i<w.
 \]
More generally, for $p,q\in \FF\langle X\rangle$, we write 
$$
 p\equiv
q\mod (S,w) 
$$ 
if $p-q\equiv0  \mod (S,w)$.

\begin{definition}
Let $S\subseteq \FF\langle X\rangle$ be monic.
If all compositions $(f,g)_w$ of $f,g\in S$ are trivial modulo $(S,w)$,
then $S$ is called a {\it Gr\"{o}bner-Shirshov basis} (with respect to $<$) in $\FF\langle X\rangle$,
or a \emph{\gsb} for the ideal $\Id(S)$,
or a \emph{\gsb} for the quotient algebra
$\FF\langle X\rangle/\Id(S)=\FF\langle X|S\rangle$.

\end{definition}

The following Composition-Diamond Lemma was first proved by Shirshov \cite{shirshov1962Lie} for Lie algebras.
Bokut \cite{bokut1976embeddings} specialized Shirshov's approach to associative
algebras (see also Bergman \cite{bergman1978diamond}) and obtain the following lemma.
For the case of commutative algebras, this lemma is known as Buchberger's Theorem
\cite{Buchberger1965}.

\begin{lemma}\label{lemma-CD}
 {\bf (Composition-Diamond Lemma)}
Let $S\subseteq \FF\langle X\rangle$ be a monic set
and let $<$ be a monomial order on $X^*$.
Then the following statements are equivalent.
\begin{enumerate}[(i)]
\item
 $S $ is a Gr\"{o}bner-Shirshov basis for $\Id(S)$ with
respect to $<$.
\item
For all $0\neq f\in \FF\langle X\rangle$,
$f\in \Id (S)\implies \overline{f}=a\overline{s}b$
for some $s\in S$ and $a,b\in  X^*$.
\item
$\Irr (S):= \{ u \in X^* \mid  u \neq a\overline{s}b ,s\in S,a ,b \in X^*\}$
is an $\FF$-basis of the algebra $\FF\langle X | S \rangle.$
\end{enumerate}
\end{lemma}

Given $u\in X^*$ and $R\subseteq \FF\langle X\rangle$,
we say $u$ is \emph{normal} with respect to $R$, or $u$ is \emph{$R$-normal}, if $u\in \Irr(R)$.
A \gsb\ $S$ is called \emph{reduced}
if for all $f\in S$ and all monomials $u$ with nonzero coefficient in $f$, we have that
$u$ is normal with respect to $S\setminus\{f\}$.
By \cite[Theorem 3]{bokut2014grobner}, every ideal of $\FF\langle X\rangle$ has a unique reduced \gsb\ for a fixed monomial order.
In applications of \gsb\ in this article, we usually use a reduced \gsb. 

Without loss of generality,
we always suppose that a reduced \gsb\ $R$ does not contain a polynomial $f$ satisfying that $\overline{f}\in X$ and $f-\overline{f}\in \FF\langle X\setminus \{\overline{f}\}\rangle$.
This assumption implies that $X\subseteq \Irr(R)$.

Let $G=\sgp XS$ be a semigroup presented by generators $X$ and defining relations $S=\{u_i=v_i\mid i\in I\}$ for some index set $I$.
We will identify a semigroup relation $u=v\ ( u,v\in X^*)$ with the algebra relation $u-v=0$ and with the polynomial (binomial) $u-v\in \FF\langle X \rangle$.
Then the semigroup algebra $\FF G$ has the
presentation $\FF G=\FF\langle X \mid S\rangle$.
Since each composition of binomials is again a binomial,
it follows from Shirshov's algorithm
that there exists a Gr\"{o}bner-Shirshov basis $R$ for $\FF G$
consisting of binomials such that $G=\sgp XR$.
Also, $R$ does not depend on the field $\FF$.
We also call $R$ a Gr\"{o}bner-Shirshov basis for the semigroup $G$.

\begin{lemma}
  \label{lemma-NF-group}
  Let $G=\sgp XS$ be a semigroup and  let $R$ be a \gsb\ for $G$.
  Then $\Irr(R)$ is a normal form of $G$.
\end{lemma}

\section{Growth of associated monomial algebras}
\label{sec-GK-MonomialAlg}
Throughout this section,
let $A=\FF\langle X\mid R\rangle=\FF\langle X \rangle/I$ be an algebra generated by a finite set $X$
unless otherwise stated,
where $I=\Id(R)$ is the ideal in $\FF\langle X\rangle$ generated by $R\subseteq \FF\langle X\rangle$.
In particular, if $R\subseteq X^*$ then $A$ is called a \emph{monomial algebra}.
Denote by $\G(A)$ the \emph{growth} of algebra $A$ (see \cite[page 6]{Krause-Lenagan_2000} for definition).

\begin{definition}
  \label{def-MonomialAlg}
Let $J$ be the ideal of $\FF\langle X\rangle$ generated by $\overline{I}=\{\overline{f}\mid f\in I\}$.
We call the monomial algebra
$$\widetilde{A_<}=\FF\langle X\rangle/J=\FF\langle X\mid \overline{I}\rangle$$
the \emph{associated monomial algebra of $A$ with respect to $<$}.
We Write $\widetilde{A}$ instead of $\widetilde{A_<}$ if no confusion arises.
\end{definition}  

Note that $\widetilde{A}=\FF\langle X\mid \overline{I}\rangle\neq \FF\langle X\mid \overline{R}\rangle$ in general
and that $\widetilde{A}=\FF\langle X\mid \overline{R}\rangle$ if $R$ is a \gsb\ for $A$.
More generally, we have the following lemma.

\begin{lemma}
  \label{lemma-GSB-MA}
Let $R'$ be a  \gsb\ for $A$ with respect to a monomial order $<$ on $X^*$.
Then 
$$
\widetilde{A_<}=\FF\langle X\mid \overline{I}\rangle=\FF\langle X\mid \overline{R'}\rangle.
$$
\end{lemma}
\begin{proof}
  Since $\Id(\overline{R'})\subseteq \Id(\overline{I})$,
  it suffices to prove that $\Id(\overline{I})\subseteq \Id(\overline{R'})$.
  Suppose $f\in \Id(\overline{I})$.
  Then $f=f_1\overline{g}f_2$ for some $f_1,f_2\in \FF\langle X\rangle$ and $g\in I$.
  Since $R'$ is a \gsb,
  it follows from Lemma \ref{lemma-CD} that
  $\overline{g}=a\overline{r}b$ for some $r\in R'$ and $a,b\in X^*$.
  Thus $f=f_1a\overline{r}bf_2\in \Id(\overline{R'})$.
\end{proof}

Under the assumption of Lemma \ref{lemma-GSB-MA},
it follows from Lemma \ref{lemma-CD} that 
$\Irr(R')=\Irr(\overline{R'})$ is an $\FF$-basis for both $A$ and $\widetilde{A}$ and thus
$A$ and $\widetilde{A}$ are isomorphic as vector spaces.
However, $A$ and $\widetilde{A}$ do not necessarily have the same GK-dimension.
For the relation between $\GK(A)$ and $\GK(\widetilde{A})$,
we have the following lemma.
\begin{lemma}\label{lem-<}
Let $<$ be a monomial order on $X^*$.
Then the following statements hold.
\begin{enumerate}[(i)]
  \item $\G(\widetilde{A})\leq \G(A)$ and thus
   $\GK(\widetilde{A})\leq \GK(A)$.
  \item If  $\GK(\widetilde{A})=\infty$, then $\GK(A)=\infty$.
  \item If  $\widetilde{A}$ has a free subalgebra of rank $2$, then so does $A$.
\end{enumerate} 
\end{lemma}
\begin{proof}
  (i) Let $V=\FF+\FF X$ and
  let $R'$ be a  \gsb\ for $A$ with respect to $<$.
  Then $V$ is a finite dimensional generating subspace of both $A$ and $\widetilde{A}$.
  We write $V$ as $V_A$ (or $V_{\widetilde{A}}$, respectively) to emphasis that 
  $V$ is thought of as a subspace of $A$ (or $\widetilde{A}$, respectively). 
  By definition,
  $(V_A)^n$ and  $(V_{\widetilde{A}})^n$, $n>0$, are both spanned by the set
  $$
  B_n:=\{y_1y_2\dots y_n\mid y_i\in  \{1\}\cup X, 1\leq i\leq n\}. 
  $$
  For each $y\in B_n$, if $y\not\in \Irr(R')$,
  then $y$ contains a subword $u\in \overline{R'}$ and thus $y=0$ in the algebra ${\widetilde{A}}$.
  Thus $C_n:=B_n\cap \Irr(R')$ is an $\FF$-basis of $(V_{\widetilde{A}})^n$.
  Note that $C_n$ is linearly independent in $(V_A)^n$.
  Hence $\dim (V_{\widetilde{A}})^n\leq \dim (V_A)^n$ for all $n>0$ and 
  thus $\G(\widetilde{A})\leq \G(A)$.
  Therefore
  $\GK(\widetilde{A})\leq \GK(A)$.
  
  (ii) It follows from part (i).
  
  (iii) If  $\widetilde{A}$ has a free subalgebra $B$ of rank $2$,
  then the algebra  $B$, as a subalgebra of $A$, is also free of rank $2$.
\end{proof}

The following lemma is well-known (see, for example, \cite[Remark 4.1]{bell2015growth}, or \cite[page 157]{Krause-Lenagan_2000} for the finitely presented case),
which is a special case of Theorem \ref{thm-main-Intro}.

\begin{lemma}\label{lem-MonomialAlg}
 Let $<$ be a deg-lex order on $X^*$.
  Then $
  \GK(A)=\GK(\widetilde{A}).
  $
\end{lemma}

The following example shows that Lemma \ref{lem-MonomialAlg} is not necessarily true for an arbitrary monomial order.
\begin{example}\label{exam-GKneq}
\label{exam_Ore}
\upshape
  Let $A=\FF\langle x,y\mid xy=y^2x\rangle$.
  Now we construct two different associated monomial algebras of $A$ with respect to two monomial orders on $X^*$
  and show that these two associated monomial algebras have different GK-dimensions.
  
  Let $<_1$ be the deg-lex order on $\{x,y\}^*$ with $y<_1x$.
  Then $xy<_1 y^2x$ and $R_1:=\{y^2x-xy\}$ is a \gsb\ for $A$ with respect to $<_1$.
  The associated monomial algebra of $A$ with respect to $<_1$ is
  $$
    \widetilde{A}_{<_1} =\FF\langle x,y\mid y^2x =0\rangle=\FF\Irr(y^2x),
  $$
  where
  \begin{align*}
    \Irr(y^2x) =\{y^sx^{i_0}yx^{i_1}y\cdots y x^{i_m}y^t\mid 0\leq s\leq 1, m\geq0,t\geq0,  \text{ each } i_j\geq 1\}.
  \end{align*}
  For any (noncommutative) polynomial $g(z_1,z_2)\in \FF\langle z_1,z_2\rangle$,
  we have that $g(yx, yx^2)\in \Irr(y^2x)$.
  Hence $\{yx, yx^2\}$  generates a free subalgebra of rank $2$ in $\widetilde{A}_{<_1}$.
  Thus it follows from Lemma \ref{lem-<} that 
  $A$ also has a free subalgebra of rank $2$.
  In particular, by Lemma \ref{lem-MonomialAlg},
  $$\GK(A)=\GK( \widetilde{A}_{<_1})=\infty.$$
  
  Now let $<_2$ be the \emph{reverse tower order} on $\{x,y\}^*$ with $y<_2x$,
  i.e.,
  for 
  \begin{align*}
    u&=y^{i_0}xy^{i_1}x\cdots y^{i_{n-1}}xy^{i_n}\in \{x,y\}^*, \\
    v&=y^{j_0}xy^{j_1}x\cdots y^{j_{m-1}}xy^{j_m}\in \{x,y\}^*,
  \end{align*}
  we set $u<_2 v $ if
  $$
    (n,i_n,i_{n-1},\dots,i_0) < (m,j_m,j_{m-1},\dots,j_0)\ \text{ lexicographically.}
  $$
  Then it is routine to check that $<_2$ is a monomial order on $\{x,y\}^*$ and $y^2x<_2 xy$.
  It is clear that $R_2:=\{xy-y^2x\}$ is a \gsb\ for $A$ with respect to $<_2$.
  Thus the associated monomial algebra of $A$ with respect to $<_2$ is
  $$
   \widetilde{A}_{<_2}=\FF\langle x,y\mid xy=0\rangle=\FF\{y^ix^j\mid i,j\geq 0\}.
  $$
  Thus $\GK( \widetilde{A}_{<_2})=2$.
  \hfil\qed
\end{example}

In particular,
we have the following remark from Example \ref{exam_Ore}.

\begin{remark}
\label{remark-tower}\upshape
 If $R$ is a \gsb\ of $A$ with respect to a reverse tower order (similarly to a tower order (Example \ref{exam-Tower}) on $X^*$,
 then $\GK(A)$ is not necessarily equal to $\GK(\widetilde{A})$.
\end{remark}

Now we are ready to prove Theorem \ref{thm-main-Intro},
which is a generalization of Lemma \ref{lem-MonomialAlg}.

\begin{proof}[Proof of Theorem \ref{thm-main-Intro}]
  Suppose $R'$ is a \gsb\ for $A$ with respect to $<$.
  The space $V:=\FF+\FF X$ is a finite dimensional generating subspace of both $A$ and $\widetilde{A}=\widetilde{A_<}$.
  We write $V$ as $V_A$ (or $V_{\widetilde{A}}$, respectively) to emphasis that 
  $V$ is thought of as a subspace of $A$ (or $\widetilde{A}$, respectively). 
  
  By definition, the subspace $(V_{A})^n$ for $n\geq 0$ is spanned by
  $$
  U:=
  \{
  x_{i_1}x_{i_2}\cdots x_{i_t}\mid 0\leq t\leq n, \ x_{i_j}\in X
  \}.
  $$
  We claim that
  each $u\in U$ is a linear combination (in $A$) of $R'$-normal
   words that are not greater than $u$.
  
  \emph{Proof of the claim.}
  We prove the claim by induction on $u$.
  If $u=1$, then the claim holds since $1\in\Irr(R')$.
  Suppose that $u>1$ and that the claim holds for words less than $u$.
  If $u\in \Irr(R')$ then the claim holds since $u$ itself is a desired linear combination.
  Assume $u\not\in \Irr(R')$.
  Then there exist $f\in R'$ and $a,b\in X^*$ such that $u=a\overline{f}b$.
  Suppose $$
  f=\alpha_1v_1+ \alpha_2 v_2+\dots+ \alpha_mv_m, \ v_i\in X^*, \ v_1>v_2>\dots>v_m, \ 0\neq \alpha_i\in \FF, \ 1\leq i\leq m.
  $$
  Then, in algebra $A$, we have that
  $$
  u=a\overline{f}b=afb-a(f-\overline{f})b=-a(f-\overline{f})b=-\alpha_2 av_2b-\alpha_3 av_3b-\dots-\alpha_m av_mb,
  $$
  where $av_ib<av_1b=u$ for $2\leq i\leq m$,
  since $<$ is a monomial order.
  By the induction hypothesis, $av_ib$ is a linear combination (in $A$) of $R'$-normal words that are not greater than $av_ib$ for $2\leq i\leq m$.
  Thus $u$ is a linear combination (in $A$) of $R'$-normal words that are not greater than $av_1b=u$. 
  Therefore the claim holds.
  
  It follows from the assumption on $f(x)$ that
  the inequality $av_ib<u\in U$ implies the inequality $|av_ib|\leq f(|u|)\leq  c|u|^d\leq cn^d$ for some $c>0$, $2\leq i\leq m$ and $|av_ib|\gg 0$.
  Thus it follows from the above claim that $(V_{A})^n$ is contained in $\FF\Irr(R')^{\leq cn^d}$, the space spanned by the set $\Irr(R')^{\leq cn^d}=\{v\in \Irr(R')\mid |v|\leq cn^d\}$.
  Note that 
  $
    \Irr(R')^{\leq n}
  $
  is a linear basis of $(V_{\widetilde{A}})^n$ for all $n\geq 0$.
  Hence $\FF\Irr(R')^{\leq cn^d}=(V_{\widetilde{A}})^{cn^d}$ as $\FF$-spaces.
  Thus,
  denoting
  $ d_A(n):=\dim (V_A)^n$ and 
  $ d_{\widetilde{A}}(n):=\dim (V_{\widetilde{A}})^{n}$ for $n\geq 0$,
  \begin{align}\label{eqn-V^n}
    d_A(n)\leq 
    d_{\widetilde{A}}(cn^d), \ \forall n\gg 0.
  \end{align}
  Therefore
  \begin{align*}
  \GK(A)&=\limsup_{n\to\infty}\frac{\ln d_A(n)}{\ln n}\\
        &\leq 
  \limsup_{n\to\infty}\frac{\ln d_{\widetilde{A}}(cn^d)}{\ln n}\ \ \ \
  \text{(by inequality (\ref{eqn-V^n}))}\\
  &=\limsup_{n\to\infty}\left(\frac{\ln d_{\widetilde{A}}(cn^d)}{\ln (cn^d)}\cdot\frac{\ln (cn^d)}{\ln n}\right)\\
  &=d\limsup_{n\to\infty}\frac{\ln d_{\widetilde{A}}(cn^d)}{\ln (cn^d)}\\
  &\leq d\limsup_{n\to\infty}\frac{\ln d_{\widetilde{A}}(n)}{\ln n}\\
  &= d\GK(\widetilde{A}).
  \end{align*}
  The inequality  $\GK(\widetilde{A})\leq \GK(A)$ follows from Lemma \ref{lem-<}.
  The proof is completed.  
\end{proof}

Theorem \ref{thm-main-Intro} has the following corollaries.
\begin{corollary}\label{cor-c}
  Suppose $<$ is a monomial order on $X^*$ such that
  $u<v$ implies $|u|\leq {c}|v|$ for some fixed real number $c>0$ and all $u,v\in X^*$.
  Then 
  $\G(A)=\G(\widetilde{A})$ and
  $
  \GK(A)=\GK(\widetilde{A}).
  $
\end{corollary}
\begin{proof}
  It follows from Lemma \ref{lem-<} that $\G(\widetilde{A})\leq \G(A)$.
  Setting $f(x)=cx$ and $d=1$ in Theorem \ref{thm-main-Intro}, 
  it follows from inequality (\ref{eqn-V^n}) that
  $
  \G(A)\leq\G(\widetilde{A}).
  $
  Hence $
  \G(A)=\G(\widetilde{A})
  $ and consequently $
  \GK(A)=\GK(\widetilde{A}).
  $
\end{proof}

A \emph{weighted deg-lex order} $<$ on $X^*$ is a generalization of the deg-lex order.
Let $\deg_w(x_i)=d_i\in \NN^*$ for $x_i\in X$.
For a non-empty word $u=x_{i_1}\dots x_{i_s}\in X^*$,
we call $\deg_w(u):=d_{i_1}+\dots +d_{i_s}$ the \emph{weighted degree} (or simply \emph{degree} if no confusion arises) of $u$.
Fix an order on each $X_m:=\{x_i\in X\mid \deg_w(x_i)=m\}$, $m\geq 1$. 
Then we compare non-empty words 
$$
u=y_1\cdots y_s,\ y_i\in X, 1\leq i\leq s
$$
and 
$$
v=z_1\cdots z_t,\ z_j\in X, 1\leq j\leq t
$$
as follows:
$u<v$ if either
$
\deg_w(u)<\deg_w(v)
$
or
$$
\deg_w(u)=\deg_w(v)
\text{ and }
(y_1,\dots, y_s)<(z_1,\dots,z_t)
\text{ lexicographically}.
$$ 
Set $1<u$ for all non-empty word $u$.
The resulting order is called a \emph{weighted deg-lex order} on $X^*$.
It is routine to check that a weighted deg-lex order is a monomial order.

\begin{corollary}\label{cor-deg-lex}
   Let $<$ be a {weighted} deg-lex order on $X^*$.
   Then $\G(A)=\G(\widetilde{A})$ and $\GK(A)=\GK(\widetilde{A}).$
\end{corollary}

\begin{proof}
  Suppose $X=\{x_1,x_2,\dots,x_m\}$ and  $\deg_w(x_i)=d_i\geq 1$, $1\leq i\leq m$.
  Let 
  $d_0:=\max\{d_i\mid 1\leq i\leq m\}.$
  Suppose $u,v\in X^*$ and $u<v$.
  Then $\deg_w(u)\leq \deg_w(v)$, 
  and thus
  $$
  |u|\leq \deg_w(u)\leq \deg_w(v)\leq d_0|v|.
  $$
  Now the statement follows from Corollary \ref{cor-c}.
\end{proof}

Let $A$ be an $\FF$-algebra and $\sigma$ be an $\FF$-endomorphism of $A$.
Recall that a linear map $\delta:A\to A$ is called a \emph{$\sigma$-derivation} if
$$
\delta(ab)=\sigma(a)\delta(b)+\delta(a)b,\ \forall a,b\in A.
$$ 
On the free left $A$-module $\oplus_{i=0}^{\infty} Ax^i$,
a ring structure is introduced by defining
$$
xa=\sigma(a)x+\delta(a), \ \forall a\in A.
$$
The resulting ring, denoted by $A[x; \sigma, \delta]$, is called a \emph{skew polynomial ring}
 or an \emph{Ore extension} over $A$. 
It was shown by Huh and Kim \cite{huh1996gelfand} that
the difference $\GK(A[x;\sigma,\delta])-\GK(A)$ can be an arbitrary positive natural number or even be infinite
(Zhao, Mo, and Zhang \cite{Zhao2018Gelfand} proved a similar pathological behavior for the GK-dimension of generalized Weyl algebras).
We refer interested readers to \cite[\S12.3]{Krause-Lenagan_2000} for more results on the GK-dimension of Ore extensions.

The algebra $A$ in Example \ref{exam_Ore} is an Ore extension over polynomial algebra $\FF [y]$,
that is, $A=\FF [y][x;\sigma,\delta]$,
where $\delta=0$ and 
$\sigma$ is an endomorphism of $\FF [y]$ defined by $\sigma(y)=y^2$.

More generally, we have the following corollary. 
Recall that an endomorphism $\sigma$ of an algebra $B$ is said to be
 \emph{locally algebraic}  if every finite dimensional subspace of $B$ 
 which contains the identity is contained in a finite dimensional subspace $U$ of $B$ such that $\sigma(U)\subseteq U$.

\begin{corollary}\label{cor-OreExt}
  Let $A=\FF [y][x;\sigma,\delta]$.
  Then the following statements are equivalent.
  \begin{enumerate}[(i)]
    \item $\GK(A)=2$.
    \item $\GK(A)<\infty$.
    \item $\deg(\sigma(y))\leq 1$. 
    \item $\sigma$ is locally algebraic.
  \end{enumerate}
\end{corollary} 

\begin{proof}
  (i)$\implies$(ii) It is obvious.
  
  (ii)$\implies$(iii)
  Suppose $d=\deg(\sigma(y))\geq 2$.
  It suffices to prove that $\GK({A})=\infty$.
  Note that 
  $A=\FF\langle x,y\mid xy=\sigma(y)x+\delta(y)\rangle$.
  Let $f:=\sigma(y)x+\delta(y)-xy$ and let
  $<$ be the weighted deg-lex order on $\{x,y\}^*$ with
  $\deg_w(x)=d_0:=\deg(\delta(y))$ and $\deg_w(y)=1$.
  Then $\deg_w(y^dx)=d+d_0$, $\deg_w(xy)=d_0+1$, and $\deg_w(\overline{\delta(y)})=\deg(\delta(y))\cdot \deg_w(y)=d_0$.
  Thus
  $$
  \overline{\sigma(y)x}=y^dx>xy>\overline{\delta(y)}
  $$ and consequently
  $\overline{f}=y^dx$.
  Since there is no composition between $f$ and $f$ itself,
  the set $\{f\}$ is a \gsb\ for $A$ with respect to $<$.
  Then, by Lemma \ref{lemma-CD}, $A$ has a basis $\Irr(f)$ which consists of polynomials
  $$
  y^{j_1}x^{i_1}\cdots y^{j_m}x^{i_m}y^{j_{m+1}},
  \ \ y^k
  $$
  where
  $0\leq j_1<d$, 
  $1\leq j_s<d$ for $2\leq s\leq m$,
  $j_{m+1}\geq0$,
  $m\geq 1$,
  $i_t\geq 1$ for $1\leq t\leq m$,
  and $k\geq 0$.
  Then, similarly to Example \ref{exam_Ore},
  $\{yx, yx^2\}$ generates a free subalgebra of rank $2$ in $A$.
  and thus $\GK({A})=\infty$.
  
  (iii)$\implies$(iv)
  Since $\{y^i\mid i\geq 0\}$ is a linear basis of $\FF[y]$,
  every finite dimensional subspace of $\FF [y]$ is contained in a \fd\ subspace $W_m:=\FF\{1,y,\dots,y^m\}$ for some $m\geq 1$.
  If $d:=\deg(\sigma(y))\leq 1$,
  then $\sigma(W_m)\subseteq W_m$ and thus $\sigma$ is locally algebraic.
  
  (iv)$\implies$(i)
  By  \cite[Theorem 1.2] {zhang1997note} (see also \cite[Theorem 12.3.3]{Krause-Lenagan_2000}),
  $\sigma$ is locally algebraic if and only if $\GK(A)=\GK(\FF[y])+1=2$.
\end{proof}

Corollary \ref{cor-OreExt} implies that the GK-dimension of an Ore extension over the polynomial algebra in one variable is either $2$ or $\infty$.
This dichotomy cannot be extended to the GK-dimension of an Ore extension over a polynomial algebra in more than one variable,
as shown by the following example.

\begin{example}
  \label{exam-OreExt2v}
   Let $A=\FF [y,z][x;\sigma,\delta]$ be an Ore extension with $\delta(y)=\delta(z)=0$.
   \begin{enumerate}[(i)]
     \item If $\sigma=\id$ (i.e., the identity map), then $\GK(A)=3$.
     \item If $\sigma(y)=y$ and $\sigma(z)=yz$,
   then $\GK(A)=4$.\item If $\sigma(y)=y^2$ and $\sigma(z)=z$, then $\GK(A)=\infty$.
   \end{enumerate}
\end{example}

\begin{proof}
  (i) In this case, $A=\FF[x,y,z]$ and thus $\GK(A)=3$.
  
  (ii)   
  Note that $A=\FF\langle x,y,z\mid xy-yx, xz-yzx, zy-yz\rangle$.
  Let $\Omega=\{y^i\mid i\geq0\}$,
  which is contained in the center of $A$.
  Then the localization $A\Omega^{-1}$ of $A$ is the differential difference algebra of GK-dimension $4$
  studied in \cite[Example 3.12]{ZhangZhao2013algebras}.
  Since $\Omega$ is a multiplicatively closed subset of regular elements that is contained in the center of $A$,
  it follows from \cite[Proposition 4.2]{Krause-Lenagan_2000} that
  $\GK(A)=\GK(A\Omega^{-1})=4$.
  
  (iii)
  Note that 
  $A=\FF\langle x,y,z\mid R\rangle$, where
  $R=\{y^2x-xy, zx-xz, zy-yz\}$ is a reduced \gsb\ of $A$
  with respect to the deg-lex order on $\{x,y,z\}^*$ with $x<y<z$.
  In a similar way as in Example \ref{exam_Ore},
  we get that
  $\{yx,yx^2\}$ generates a free subalgebra of rank $2$ in $A$
  and thus $\GK(A)=\infty$.  
\end{proof}
  
\section{Growth of Manturov groups}
\label{sec-Manturov}
In this section, we study the growth and GK-dimension of Manturov groups via the method of \gsbs.

\subsection{Definition of Manturov groups}
Let us first recall the definition of the group $G_n^k$ from \cite{manturov2015non,manturov2015on}.
Given two positive integers $n>k$,
Denote by $I_n^k$ the set of all sets of $k$ distinct indices from $\{1,2,\ldots,n\}$.
Clearly, the cardinality of $I_n^k$ is ${n\choose k}$.
The \emph{Manturov $(k,n)$-group} $G_n^k$ is defined by generators and relations as follows.
\begin{enumerate}[(a)]
  \item Generators: $G_n^k$ has $n\choose k$ generators $a_m$, $m\in I_n^k$.
\item Relations:
there are three types of relations for $G_n^k$.

(i) (Involution relations) $a_m^2=1$, for all $m\in I_n^k$.

(ii) (Far commutativity relations)
$a_ma_{m'}=a_{m'}a_m$, for any $m,m'\in I_n^k$ with $|m\cap m'|<k-1$.

(iii) (Tetrahedron relations) For each $(k+1)$-tuple $U$ of indices $u_1,\dots,u_{k+1}\in \{1,\dots,n\}$,
consider the $k+1$ subsets $m^j=U\setminus\{u_j\}\in I_n^k$, $j=1,\dots,k+1$.
With $U$, we associate the relation:
\begin{align*}\label{rel_tetrahedron}
  a_{m^1}a_{m^2}\cdots a_{m^{k+1}}=a_{m^{k+1}}\cdots a_{m^2}a_{m^1}.
\end{align*}
\end{enumerate}

Note that if $n=k+1$ or if $k=1$, then the far commutativity relations never happen.
In the tetrahedron relations,
$|m^i\cap m^j|=k-1$ for $i\neq j$ 
and thus the far commutativity relation does not happen to these $m^i$ and $m^j$;
for two tuples $U$ and $U'$ that differ only by order reversal, we get the same relation,
and hence the group $G_n^k$ has totally $\frac{{n\choose k+1}(k+1)!}{2}$ tetrahedron relations.
\begin{example}
  \label{exam_G1n}
  \upshape
  \begin{enumerate}[(i)]
    \item 
    The smallest possible $G_n^k$ is $G_2^1$,
  which has two generators, denoted by $a$ and $b$.
  Then
  \[
  G_2^1=\gp{ a,b}{a^2=b^2=1,ab=ba}\simeq \ZZ_2\oplus\ZZ_2.
  \]
   \item 
    More generally, for $n\geq 1$,
    \begin{align*}
      G_n^1&=\gp{a_i}{a_i^2=1, a_ta_s=a_sa_t, \ 1\leq i\leq n, 1\leq s<t\leq n}\\
      &\simeq \ZZ_2\oplus\ZZ_2\oplus\cdots\oplus\ZZ_2 \ \ (n \text{  copies of } \ZZ_2).
    \end{align*}
  \end{enumerate}
\end{example}


When $k>1$, the group $G_n^k$ seems complicated.
Even when $k=2$, the groups $G_n^2$ are ``deeply non-trivial'' \cite{manturov2017the}.
\subsection{Gr\"obner-Shirshov basis, word problem, normal form, and growth of $G_3^2$}
Let $a=a_{12}$, $b=a_{13}$, and $c=a_{23}$.
Then the group $G_3^2$ has a group presentation
\[
G_3^2=\gp{ a,b,c}{ a^2=b^2=c^2=1, bca=acb, cab=bac, cba=abc},
\]
which is also a semigroup presentation,
i.e.,
\begin{align*}
  G_3^2=\sgp{ a,b,c}{ a^2=b^2=c^2=1, bca=acb, cab=bac, cba=abc}.
\end{align*}

Bai and Chen \cite[Theorem 3]{BC17} (see also \cite[Theorem 1]{AlHussein2021}) obtained a finite \gsb\ for $G^2_3$ with respect to a tower order.
However, it follows from Remark \ref{remark-tower} that
a \gsb\ with respect to a tower order cannot be used directly 
to calculate the GK-dimension of $\FF G^2_3$.
The following theorem gives an infinite \gsb\ for $G^2_3$ with respect to a deg-lex order. 

\begin{theorem}
  \label{thm_gsb-G32}
  Let $<$ be the deg-lex order on $\{a,b,c\}^*$ such that $a<b<c$.
  Then the set $R'$ consisting of the following  relations is a \gsb\ for $G_3^2$ with respect to $<$:
   \begin{enumerate}[(i)]
     \item $x^2-1,\ x\in\{a,b,c\}$;
     \item $bca-acb;$
     \item $cab-bac$;
     \item $cba-abc;$
     \item $ b(ac)^mb-(ca)^m$, $ m\geq 1$.
   \end{enumerate}
\end{theorem}
\begin{proof}
By an induction on $m$,
we can show that relations (v) in $R'$ can be deduced from relations (i)-(iv).
Hence $ G_3^2=\sgp{a,b,c}{ R'}$.
Now it suffices to prove that all compositions of polynomials in $R'$ are trivial.

Note that there is no inclusion composition in $R'$.
The intersection compositions in $R'$ and all ambiguities $w$ are list as follows (where we use, for example, (ii)$\wedge $(iii) to denote the intersection composition(s) of relations (ii) and (iii)):
\begin{itemize}
  \item
  (i)$\wedge $(i):
  $w=xxx$, \ $x=a,b,c$;  \\
  (i)$\wedge $(ii): $w=bbca$;\\
  (i)$\wedge $(iii): $w=ccab$;\\
  (i)$\wedge $(iv): $w=ccba$;\\
  (i)$\wedge $(v): $w=bb(ac)^mb$, $m\geq 1$;
  \item
 (ii)$\wedge $(i): $w=bcaa$;\\
  (ii)$\wedge $(iii): $w=bcab$;
  \item
  (iii)$\wedge $(i): $w=cabb$;\\
  (iii)$\wedge $(ii): $w=cabca$;\\
  (iii)$\wedge $(v): $w=cab(ac)^mb$, $m\geq 1$;
  \item
  (iv)$\wedge $(i): $w=cbaa$;\\
  (iv)$\wedge $(v): $w=cbac(ac)^{m-1}b$, $m\geq 1$;
  \item
  (v)$\wedge $(i): $w=b(ac)^mbb$, $m\geq 1$;\\
  (v)$\wedge $(ii): $w=b(ac)^mbca$, $m\geq 1$;\\
  (v)$\wedge $(iv): $w=b(ac)^{m-1}acba$, $m\geq 1$;\\
  (v)$\wedge $(v): $w=b(ac)^mb(ac)^lb$, \ $m,l\geq 1$.
  \end{itemize}

It is routine to check that each of above compositions is trivial modulo $(R',w)$.
We list several typical ones here,
the other compositions can be checked similarly.
\begin{align*}
  \text{(i)}\wedge \text{(v)}:\ \ \ \ \ &(b^2-1)(ac)^mb-b(b(ac)^mb-(ca)^m)\\
=&-(ac)^mb+b(ca)^m\\
\equiv&-(ac)^mb+(ac)^mb\ \ \ \ \ \ \ \mbox{(using (ii) $m$ times)}\\
=&0.
\end{align*}
\begin{align*}
  \text{(v)}\wedge \text{(iv)}:\ \ \ \ \ &(b(ac)^mb-(ca)^m)a-b(ac)^{m-1}a(cba-abc)\\
=&-(ca)^ma+b(ac)^{m-1}aabc\\
\equiv&-(ca)^{m-1}c+b(ac)^{m-1}bc\ \ \ \ \ \ \ \mbox{(using (i))}\\
\equiv&-(ca)^{m-1}c+(ca)^{m-1}c\ \ \ \ \ \ \ \mbox{(using (v))}\\
=&0.
\end{align*}
\begin{align*}
  \text{(v)}\wedge \text{(v)}:\ \ \ \ \ h:=&(b(ac)^mb-(ca)^m)(ac)^lb-b(ac)^{m}(b(ac)^lb-(ca)^l)\\
= &-(ca)^m(ac)^lb+b(ac)^{m}(ca)^l.
\end{align*}
If $m=l$, then $h=0$.
If $m>l$, then
\begin{align*}
h\equiv&-(ca)^{m-l}b+b(ac)^{m-l}\\
\equiv&-b(ac)^{m-l}+b(ac)^{m-l}\ \ \ \ \ \ \ \mbox{(using (iii)  $m-l$ times)}\\
=&0.
\end{align*}
If $m<l$, then
\begin{align*}
h\equiv&-(ac)^{l-m}b+b(ca)^{l-m}\\
\equiv&-(ac)^{l-m}b+(ac)^{l-m}b\ \ \ \ \ \ \ \mbox{(using (ii) $l-m$ times)}\\
=&0.
\end{align*}
Therefore, $R'$ is a \gsb\ of the group $G^2_3$.
\end{proof}

Note that $R'$ contains infinitely many relations.
However, given a word $u\in\{a,b,c\}^*$,
there are only finitely many relations in $R'$ 
that can be used to reduce $u$.
Moreover, if $u\not\in \Irr(R')$,
then an element $f=\overline{f}-f'\in R'$ such that $u=a\overline{f}b$ can be determined, where $f',a,b\in X^*$.
We can do the reduction 
$u\longrightarrow u_1$, where $u_1=af'b<u$ as words and $u_1=u$ as elements in $G^2_3$.
If $u_1\not\in \Irr(R')$, then we can continue the reduction $u_1\longrightarrow u_2$ such that $u_2<u_1$ as words and $u_2=u_1$ as elements in $G^2_3$.
The chain of reductions of $u$ will terminate in finitely many steps since $<$ is a well order, say
$
u\longrightarrow u_1\longrightarrow u_2\longrightarrow \dots \longrightarrow u_k,
$
where $k\in \NN$ and $u_k\in \Irr(R')$.
Similarly, for $v\in X^*$, we have 
$
v\longrightarrow v_1\longrightarrow v_2\longrightarrow \dots \longrightarrow v_l,
$
where $l\in \NN$ and $u_l\in \Irr(R')$.
Since $\Irr(R')$ is a normal form of the group $G^2_3$,
we obtain that $u=v$ in $G^2_3$ if and only if $u_k=v_l$ as words.
That is, we have the following corollary
(which also follows from the \gsb\ obtained by Bai and Chen \cite{BC17}).

\begin{corollary}\label{cor-WordProb23}
  The word problem for the group $G_3^2$ is solvable.
\end{corollary}

Write $u\sqsubseteq v$ if $u$ is a subword of $v$.
The following lemma gives a normal form of the group $G^2_3$.
\begin{corollary}\label{coro-NormalForm}
  The group $G^2_3$ has a normal form
  \begin{align*}
    N:=\{u\mid u\sqsubseteq (ab)^s(ac)^t,(ab)^s(cb)^t \text{ or } (ac)^s(bc)^t, \ \ s,t\in\NN\}.
  \end{align*}
\end{corollary}
\begin{proof}
  By Theorem \ref{thm_gsb-G32}, $R'$ is a \gsb\ of $G^2_3$.
  Thus, by Lemma \ref{lemma-NF-group}, 
  $\Irr(R')$ is a normal form of $G^2_3$.
  Hence it suffices to prove that $N=\Irr(R')$.
  
  Note that $(ab)^s(ac)^t, (ab)^s(cb)^t,(ac)^s(bc)^t\in \Irr(R')$ for all $s,t\in\NN$.
  Thus all subwords of $(ab)^s(ac)^t, (ab)^s(cb)^t,(ac)^s(bc)^t$ belong to $\Irr(R')$ for $s,t\in\NN$.
  Hence $N\subseteq \Irr(R')$.  
    
  Now suppose $u\in\Irr(R')$. 
  Next we prove $u\in N$ by induction on $|u|$.
  If $|u|=1,2$, then $u\in\{a,b,c,ab,ac,ba,bc,ca,cb\}\subseteq N$.
  Suppose for all $u$ with $|u|=n\geq 2$ we have $u\in N$.
  We consider the case $|u|=n+1$.
  Write $u=vx$, where $|v|=n\geq 2$ and $x\in\{a,b,c\}$.
  By the induction hypothesis, $v\in N$.
  According to the first letter of $v$,
  we have three cases to consider.
    
  (Case 1) Suppose the first letter of $v$ is $a$,
  i.e., $v=a\cdots$.
  
  (Case 1.1)
   Suppose $v=ab\cdots$.
   Then $v\sqsubseteq (ab)^s(ac)^t$ or $v\sqsubseteq (ab)^s(cb)^t$ for some $s,t\in \NN$ 
   and thus $v$ is one of the following words:
   \begin{align*}
     v &= (ab)^i\sqsubseteq (ab)^s,\ 1\leq i\leq s,\\
     v &= (ab)^ia\sqsubseteq (ab)^s(ac)^t,\ 1\leq i\leq s,\\
     v &= (ab)^i(ac)^j\sqsubseteq (ab)^s(ac)^t,\ 1\leq i\leq s, 1\leq j\leq t,\\
     v &= (ab)^i(ac)^ja\sqsubseteq (ab)^s(ac)^t,\ 1\leq i\leq s, 1\leq j<t,\\
     v &= (ab)^i(cb)^jc\sqsubseteq (ab)^s(cb)^t,\ 1\leq i\leq s, 0\leq j<t,\\
     v &= (ab)^i(cb)^j\sqsubseteq (ab)^s(cb)^t,\ 1\leq i\leq s, 1\leq j\leq t.
   \end{align*}
  It is routine to check that $u=vx\in N$ for each possibility of $v$.
  For instance,
  if $v=(ab)^i$,
  then, noting that $vb\not\in\Irr(R')$, we have that either $u=vx=va=(ab)^ia\in N$ or
  $u=vx=vc=(ab)^ic\in N$;
  if $v=(ab)^i(cb)^j$,
  then $va,vb\not\in\Irr(R')$ and thus $u=vx=vc=(ab)^i(cb)^jc\in N$.
  
  (Case 1.2)
  Suppose $v=ac\cdots$.
  Then $v\sqsubseteq (ac)^s(bc)^t$ for some $s,t\in \NN$.
  Thus $v$ is one of the following words:
   \begin{align*}
     v &= (ac)^i\sqsubseteq (ac)^s,\ 1\leq i\leq s,\\
     v &= (ac)^ia\sqsubseteq (ac)^s,\ 1\leq i<s,\\
     v &= (ac)^i(bc)^j\sqsubseteq (ac)^s(bc)^t,\ 1\leq i\leq s, 1\leq j\leq t,\\
     v &= (ac)^i(bc)^jb\sqsubseteq (ac)^s(bc)^t,\ 1\leq i\leq s, 0\leq j<t.
   \end{align*}
  It is routine to check that $u=vx\in N$ for each possibility of $v$.

  (Case 2) Suppose $v=b\cdots$.
  
  (Case 2.1)
   Suppose $v=ba\cdots$.
   Then $v\sqsubseteq (ab)^s(ac)^t$ or $v\sqsubseteq (ab)^s(cb)^t$ for some $s,t\in \NN$.
   If $|v|=2$ then $v=ba$ and $u=bab,bac\in N$.
   Suppose $|v|>2$. 
   Then $v=bv'$, where $v'=a\cdots$ and $|v'|\geq 2$.
   Thus, by Case 1.1, $v'$ is one of the following words:
   \begin{align*}
     v' &= (ab)^i\sqsubseteq (ab)^s,\ 1\leq i\leq s,\\
     v' &= (ab)^ia\sqsubseteq (ab)^s(ac)^t,\ 1\leq i\leq s,\\
     v' &= (ab)^i(ac)^j\sqsubseteq (ab)^s(ac)^t,\ 1\leq i\leq s, 1\leq j\leq t,\\
     v' &= (ab)^i(ac)^ja\sqsubseteq (ab)^s(ac)^t,\ 1\leq i\leq s, 1\leq j<t,\\
     v' &= (ab)^i(cb)^jc\sqsubseteq (ab)^s(cb)^t,\ 1\leq i\leq s, 0\leq j<t,\\
     v' &= (ab)^i(cb)^j\sqsubseteq (ab)^s(cb)^t,\ 1\leq i\leq s, 1\leq j\leq t.
   \end{align*}
  Similarly,
  one can check that $u=bv'x\in N$ for each possibility of $v'$.
  
  (Case 2.2)
  Suppose $v=bc\cdots$.
  Then $v\sqsubseteq  (ac)^s(bc)^t$ for some $s,t\in \NN$.
  Thus $v$ is one of the following words:
   \begin{align*}
     v &= (bc)^i\sqsubseteq (ac)^s(bc)^t,\ 1\leq i\leq t,\\
     v &= (bc)^ib\sqsubseteq (ac)^s(bc)^t,\ 1\leq i<t.
   \end{align*}
  It is routine to check that $u=vx\in N$ for each possibility of $v$.

  (Case 3) Suppose $v=c\cdots$.
  Similarly to Cases 1 and 2, we obtain that
  if $v=ca\cdots$, then $v$ is one of the following words:
   \begin{align*}
     v &= (ca)^i\sqsubseteq (ac)^s(bc)^t,\ 1\leq i< s,\\
     v &= (ca)^ic(bc)^j\sqsubseteq (ac)^s(bc)^t,\ 1\leq i<s, 0\leq j\leq t,\\
     v &= (ca)^ic(bc)^jb\sqsubseteq (ac)^s(bc)^t,\ 1\leq i<s, 0\leq j<t;
   \end{align*}
  and that if $v=cb\cdots$, then $v$ is one of the following words:
  \begin{align*}
     v &= (cb)^i\sqsubseteq (ac)^s(bc)^t,\ 1\leq i\leq t,\\
     v &= (cb)^ic \sqsubseteq (ac)^s(bc)^t,\ 1\leq i\leq t.
  \end{align*}
  It is routine to check that $u=vx\in N$ for each possibility of $v$. 
  This completes the proof.
\end{proof}

\begin{corollary}\label{cor-G23}
\begin{enumerate}[(i)]
  \item The group $G_3^2$ has growth equal to $2$ (i.e. quadratic growth). In particular,
   $\GK(\FF G_3^2)=2$.
  \item The group $G^2_3$ is a finite extension of $\ZZ^2$.
  \item The group algebra $\FF G^2_3$ is a finite module over the Laurent series ring $\FF[x^{\pm1},y^{\pm1}]$.
\end{enumerate}
\end{corollary}
\begin{proof}
  (i)
  Let $A:=\FF G^2_3=\FF\langle X\mid R'\rangle$,
  where $X=\{a,b,c\}$ and $R'$ is the \gsb\ (with respect to a deg-lex order) defined in Theorem \ref{thm_gsb-G32}.
  It follows from Corollary \ref{cor-deg-lex} that 
  $$\GK(A)=\GK(\widetilde{A})
 =\limsup_{n\to\infty}\log_n \sum_{i=0}^nf(i),$$
  where
  $f(i):=\#\{u\mid u\in\Irr(R'), |u|=i\}$ for $ 0\leq i\leq n.$

 Suppose $u\in\Irr(R')$ and $|u|=2n>0$. By Corollary \ref{coro-NormalForm}, $u$ is a subword of $(ab)^s(ac)^t$, $(ab)^s(cb)^t$ or $(ac)^s(bc)^t$ for some $s,t\in \NN$.
 Suppose $u\sqsubseteq (ab)^s(ac)^t$.
 Then all distinct possibilities of $u$ are as follows:
 $$
(ba)^n, (ab)^n, (ab)^{n-1}ac,b(ab)^{n-2}(ac)a,\dots,
 (ab)(ac)^{n-1},b(ac)^{n-1}a,(ac)^{n}, (ca)^n.
 $$
 Thus $\#\{u\mid u\sqsubseteq (ab)^s(ac)^t, |u|=2n\}=2n+3$.
  Similarly, we have 
  $$
  \#\{u\mid u\sqsubseteq (ab)^s(cb)^t, |u|=2n\}=\#\{u\mid u\sqsubseteq (ac)^s(bc)^t, |u|=2n\}=2n+2.
  $$
  Thus, for $n>0$,
  \begin{align*}
   f(2n) &\geq \# \{u\mid u\sqsubseteq (ab)^s(ac)^t, |u|=2n\}= 2n+3,
 \end{align*}
 and
 \begin{align*}
   f(2n) & =\# \{u\mid u\in\Irr(R'), |u|=2n\}\\
   & \leq \# \{u\mid u\sqsubseteq (ab)^s(ac)^t, |u|=2n\}+\# \{u\mid u\sqsubseteq (ab)^s(cb)^t, |u|=2n\}\\
   &\ \ \ +\# \{u\mid u\sqsubseteq (ac)^s(bc)^t, |u|=2n\}\\
    &= 6n+7.
 \end{align*}
 Thus $f(n)\sim n$ (see \cite[page 5]{Krause-Lenagan_2000} for the definition of the equivalence relation $\sim$).
 Therefore
 $$\GK(\FF G_3^2)
 =\limsup_{n\to\infty}\log_n \sum_{i=0}^nf(i)=2$$
 and $G_3^2$ has growth equal to $2$.
 
(ii)
According to Gromov's theorem, a group of polynomial growth is a virtually nilpotent group, 
i.e., it has a nilpotent subgroup of finite index.
The growth degree of a nilpotent group is determined by the Bass-Guivarc'h formula
\cite[Theorem 11.14]{Krause-Lenagan_2000}, 
which depends on the torsion free ranks of abelian groups occurring in the lower central series. 
Then it is easy to see that a group with quadratic growth is 
in fact virtually abelian of rank $2$,
and so is a finite extension of $\ZZ^2$ (see \cite[page 196]{Krause-Lenagan_2000}).
 
(iii)
It follows from (ii).
\end{proof}

\subsection{Growth of $G_4^3$}
By definition, 
the Manturov group $G_4^3$ has the semigroup presentation
\[
G_4^3=\sgp{a,b,c,d}{R},
\]
where $ a = a_{123}, b = a_{124}, c =a_{134}, d = a_{234}$
and $R$ consists of the following relations:
\begin{enumerate}[(i)]
  \item The involution relations: $x^2=1$, $x\in\{a,b,c,d\}$.
  \item The tetrahedron relations:
  $xyzw=wzyx$, where $\{x,y,z,w\}=\{a,b,c,d\}$.
\end{enumerate}

\begin{proposition}\label{pro-Growth-G34}
  \begin{enumerate}[(i)]
  \item The group $G^3_4$ has a free subgroup of rank $2$.
  \item $\GK(\FF G^3_4)=\infty.$ 
  \item The group $G^k_n$ has a free subgroup of rank $2$ for all $n>k\geq 3$.
  \item $\GK(\FF G^k_n)=\infty$ for all $n>k\geq 3$.
\end{enumerate}  
\end{proposition}

\begin{proof}
  (i)
  It follows from Bai and Chen \cite[Theorem 3]{BC17} (see also \cite[Theorem 1]{AlHussein2021}) that $G^3_4$ has a \gsb\footnote{Note that our notation is different from that in \cite{BC17}, e.g., $a=a_{234}$ in \cite{BC17} while $a=a_{123}$ in our paper.
   However, this difference does not affect the presentation of $G^3_4$ due to the symmetry of the generators.}
   \begin{align*}
     S=& \{a^2-1, b^2-1, c^2-1, d^2-1,\\
     &ba-cdabcd, bcda-adcb, bca- dacbd, bda-cadbc, dca-cdabcdcdb, dcda-cabdcdcb\}
   \end{align*}
   with respect to the tower order on $\{a,b,c,d\}^*$ with $a>b>c>d$.
   Thus
  $$
  \Irr(S)=\{
  u_0u_1\cdots u_nv\mid
  n\geq 0, u_0\in B_1, u_j\in B_2, 1\leq j\leq n, v\in B_3
  \}
  $$
  where
  \begin{align*}
    B_1 & = \{1, a, ca, da, cda\}, \\
    B_2 & = \{ca, da, cda\}, \\
    B_3   & = \{x_1x_2\cdots x_t\mid t\geq 0, x_i\in \{b,c,d\}, 1\leq i\leq t, x_j\neq x_{j+1}, 1\leq j\leq t-1\} .
  \end{align*}
   For any polynomial $g(z_1,z_2)\in \FF\langle z_1,z_2\rangle$,
  we have that $g(ca, da)\in \Irr(S)$.
  Thus $\{ca, da\}$  generates a free subalgebra (free subgroup, respectively) of rank $2$ in $\FF G^3_4$ (in $G^3_4$, respectively).
    
  (iii)
 It follows from \cite{manturov2015non} that,
 for each $l = 1, \dots, n$,
 there is a group homomorphism (called \emph{index forgetting homomorphism}) 
 $f_l$ from $G^k_n$ onto $G^{k-1}_{n-1}$.
 By part (i), $G^3_4$ has a free subgroup of rank $2$.
 Thus $G^4_5$ has a free subgroup of rank $2$.
 By an induction on $n\geq 4$, we obtain that $G^{n-1}_n$ has a free subgroup of rank $2$ for all $n\geq 4$.
 Similarly,
 the existence of the \emph{strand-deletion homomorphism} \cite{manturov2015non}
  $d_n$ from $G^k_n$ onto $G^k_{n-1}$
  implies that $G^3_n$ has a free subgroup of rank $2$ for all $n\geq 4$.
  Combining the two types of homomorphisms,
  we obtain that
  $G^k_n$ has a free subgroup of rank $2$ for all $n>k\geq 3$.
  
  Parts (ii) and (iv) follow directly from (i) and (iii) respectively.
\end{proof}

Now we are ready to prove Theorem \ref{thm-main2}.

\begin{proof}[Proof of Theorem \ref{thm-main2}]
  (i) By Example \ref{exam_G1n}, $G^1_n$ is a finite group and thus
   $\GK(\FF G^1_n)=0$  for all  $n>1$.
   That is, $G^1_n$ has growth equal to $0$ for $n>1$.
  
  (ii) It follows form Corollary \ref{cor-G23}.
  
  (iii) It follows form Proposition \ref{pro-Growth-G34}.  
\end{proof}

\subsection{Growth of group $G^2_n$}
We have obtained the growth and GK-dimension of all Manturov groups $G^k_n$ for $n>k\geq 1$ except $G^2_n$ for $n\geq 4$.
The group $G^2_4$ has $6$ generators:
$$
a=a_{12},
b=a_{13},
c=a_{14},
d=a_{23},
e=a_{24},
f=a_{34}.
$$
By definition, 
\begin{align*}
  G^2_4=\gp{a,b,c,d,e,f}{R},
\end{align*}
where $R$ consists of the following relations
\begin{enumerate}[(i)]
  \item $a^2=b^2=c^2=d^2=e^2=f^2=1$,
  \item $fa=af, eb=be, dc=cd$,
  \item $dba=abd,dab=bad,bda=adb;
   eca=ace, eac=cae, cea=aec;\\
   fcb=bcf,fbc=cbf,cfb=bfc;
   fed=def,fde=edf,efd=dfe.$
\end{enumerate}
As far as we know,
it is quite difficult to find a reduced \gsb\ for $G^2_4$.
At present we do not know the growth of $G_4^2$.
We propose the following conjecture.

\begin{conjecture}\label{conj}
  The group $G^2_4$ has a free subgroup of rank $2$.
\end{conjecture}

Based on observations, 
it seems that $\{ab, ac\}$ generates a free subgroup of rank $2$ of $G^2_4$.
If Conjecture \ref{conj} is valid,
then it follows from the strand-deletion homomorphism (see the proof for Proposition \ref{pro-Growth-G34} (iii)) that
 $G^2_n$ has a free subgroup of rank $2$ for all $n\geq 4$.

\subsection*{Acknowledgments}
The author thanks the anonymous referees for their careful reading and valuable comments, which have improved the readability of this article.
Part of the work is motivated by a series of talks given by L.A. Bokut in South China Normal University (2016) and a conference talk given by V.O. Manturov in ICAC 2017 in Hong Kong.
We express our gratitude to Bokut and Manturov
 for fruitful and stimulation discussions.
The author would like to thank Caiheng Li and the Department of Mathematics at Southern University of Science and Technology for the hospitality during his visits.
The author wishes to thank Yuxiu Bai, Yu Li, Wenchao
Zhang, Zerui Zhang, and Fangting Zheng for helpful discussions.
 
\subsection*{Funding}
The author was partially supported by the Guangdong Basic and Applied Basic Research Foundation (Nos. 2023A1515011690, 2022A1515110634) and
the Characteristic Innovation Project of Guangdong Provincial Department of Education (No. 2023KTSCX145).

\end{document}